\documentclass[12pt]{article}
\usepackage{amssymb,latexsym,amsmath}

\begin{document}
\pagestyle{plain} \headheight=5mm \topmargin=-5mm 

\title{Birational invariants defined by Lawson homology   }
\author{ Wenchuan Hu
 }

\maketitle
\newtheorem{Def}{Definition}[section]
\newtheorem{Th}{Theorem}[section]
\newtheorem{Prop}{Proposition}[section]
\newtheorem{Not}{Notation}[section]
\newtheorem{Lemma}{Lemma}[section]
\newtheorem{Rem}{Remark}[section]
\newtheorem{Cor}{Corollary}[section]

\def\s{\section}
\def\ss{\subsection}

\def\d{\begin{Def}}
\def\t{\begin{Th}}
\def\p{\begin{Prop}}
\def\n{\begin{Not}}
\def\la{\begin{Lemma}}
\def\r{\begin{Rem}}
\def\c{\begin{Cor}}
\def\ee{\begin{equation}}
\def\aa{\begin{eqnarray}}
\def\ya{\begin{eqnarray*}}
\def\bd{\begin{description}}

\def\ed{\end{Def}}
\def\et{\end{Th}}
\def\epo{\end{Prop}}
\def\en{\end{Not}}
\def\el{\end{Lemma}}
\def\er{\end{Rem}}
\def\ec{\end{Cor}}
\def\eee{\end{equation}}
\def\eaa{\end{eqnarray}}
\def\ey{\end{eqnarray*}}
\def\ebd{\end{description}}

\def\nn{\nonumber}
\def\bp{{\bf Proof.}\hspace{2mm}}
\def\qe{\hfill$\Box$}
\def\lj{\langle}
\def\rj{\rangle}
\def\dd{\diamond}
\def\ox{\mbox{}}
\def\lb{\label}
\def\rel{\;{\rm rel.}\;}
\def\vp{\varepsilon}
\def\ep{\epsilon}
\def\mod{\;{\rm mod}\;}
\def\exp{{\rm exp}\;}
\def\Lie{{\rm Lie}}
\def\dim{{\rm dim}}
\def\im{{\rm im}\;}
\def\Lag{{\rm Lag}}
\def\Gr{{\rm Gr}}
\def\span{{\rm span}}
\def\Spin{{\rm Spin}}
\def\sign{{\rm sign}\;}
\def\Supp{{\rm Supp}\;}
\def\Sp{{\rm Sp}\;}
\def\ind{{\rm ind}\;}
\def\rank{{\rm rank}\;}
\def\Sg{{\Sp(2n,\C)}}
\def\Na{{\cal N}}
\def\det{{\rm det}\;}
\def\dist{{\rm dist}}
\def\deg{{\rm deg}}
\def\Tr{{\rm Tr}\;}
\def\ker{{\rm ker}\;}
\def\Vect{{\rm Vect}}
\def\H{{\bf H}}
\def\K{{\rm K}}
\def\R{{\bf R}}
\def\C{{\bf C}}
\def\Z{{\bf Z}}
\def\N{{\bf N}}
\def\F{{\bf F}}
\def\Da{{\bf D}}
\def\A{{\bf A}}
\def\La{{\bf L}}
\def\x{{\bf x}}
\def\y{{\bf y}}
\def\Ga{{\cal G}}
\def\Ha{{\cal H}}
\def\L{{\cal L}}
\def\Pa{{\cal P}}
\def\Ua{{\cal U}}
\def\E{{\rm E}}
\def\J{{\cal J}}

\def\m{{\rm m}}
\def\ch{{\rm ch}}
\def\gl{{\rm gl}}
\def\Gl{{\rm Gl}}
\def\Sp{{\rm Sp}}
\def\sf{{\rm sf}}
\def\U{{\rm U}}
\def\O{{\rm O}}
\def\F{{\rm F}}
\def\P{{\rm P}}
\def\D{{\rm D}}
\def\T{{\rm T}}
\def\Sa{{\rm S}}

\begin{center}{\bf \tableofcontents}\end {center}

\begin{center}{\bf Abstract }\end {center}

New birational invariants for a projective manifold are defined by
using Lawson homology. These invariants can be highly nontrivial
even for projective threefolds. Our techniques involve the weak
factorization theorem of Wlodarczyk and tools developed by
Friedlander, Lawson, Lima-Filho and others. A blowup formula for
Lawson homology is given in a separate section. As an application,
we show that for each $n\geq 5$, there  is a smooth  rational
variety $X$ of dimension $n$ such that the Griffiths groups ${\rm
Griff}_p(X)$ are infinitely generated even modulo torsion for all
$p$ with $2\leq p\leq n-3$.

\s {Introduction}

\hskip .2in In this paper, all varieties are defined over
$\mathbb{C}$. Let $X$ be an n-dimensional projective variety. The
\textbf{Lawson homology} $L_pH_k(X)$ of $p$-cycles is defined by
$$L_pH_k(X) := \pi_{k-2p}({\cal Z}_p(X)) \quad for\quad k\geq 2p\geq 0,$$
where ${\cal Z}_p(X)$ is provided with a natural topology (cf.
\cite{Friedlander1}, \cite{Lawson1}). For general background, the
reader is referred to \cite{Lawson2}.

In \cite{Friedlander-Mazur}, Friedlander and Mazur showed that there
are  natural transformations, called \textbf{cycle class maps}
 $$ \Phi_{p,k}:L_pH_{k}(X)\rightarrow H_{k}(X). $$
 Define $$L_pH_{k}(X)_{hom}:={\rm ker}\{\Phi_{p,k}:L_pH_{k}(X)\rightarrow
 H_{k}(X)\}.$$

The \textbf{Griffiths group} of codimension $q$-cycles is defined to
$$
{\rm Griff}^q(X):={\cal Z}^q(X)_{hom}/{\cal Z}^q(X)_{alg}
$$

It was proved by Friedlander \cite{Friedlander1} that, for any
smooth projective variety $X$, $L_pH_{2p}(X)\cong {\cal
Z}_p(X)/{\cal Z}_p(X)_{alg}$. Therefore

\begin{eqnarray*}
L_pH_{2p}(X)_{hom}\cong {\rm Griff}_p(X),
\end{eqnarray*} where ${\rm
Griff}_p(X):={\rm Griff}^{n-p}(X)$.

\medskip

 The main result in this paper is the following
 {\Th If $X$ is a smooth n-dimensional projective variety, then
 $L_{1}H_{k}(X)_{hom}$ and
  $L_{n-2}H_{k}(X)_{hom}$ are  smooth birational invariants for $X$. More
  precisely, if $\varphi:X \rightarrow X^{\prime}$ is a birational map
  between smooth projective manifolds $X$ and $X^{\prime}$, then
   $\varphi$ induces isomorphisms $L_{1}H_{k}(X)_{hom}\cong L_{1}H_{k}(X^{\prime})_{hom}$ for
  $k\geq 2$ and  $L_{n-2}H_{k}(X)_{hom}\cong L_{n-2}H_{k}(X^{\prime})_{hom}$ for
  $k\geq 2(n-2)$. In particular, $L_{1}H_{k}(X)_{hom}=0$  and $L_{n-2}H_{k}(X)_{hom}=0$ for any
  smooth rational variety.
  }

{\Cor Let $X$ be a smooth rational projective variety with
$\dim(X)\leq 4$, then $\Phi_{p,k}: L_pH_k(X)\rightarrow H_k(X)$ is
injective for all $k\geq 2p\geq 0$.}

{\Rem In general, for $2\leq p\leq n-3$, $L_{p}H_{k}(X)_{hom}$ is
\emph{not} a birational invariant for the smooth projective variety
$X$. This follows from the blowup formula in Lawson homology (See
Corollary 1.2, 1.3).}

{\Rem If $p=0,n-1, n$, then $L_pH_k(X)_{hom}=0$ for all $k\geq 2p$.
In these cases, the statement in the theorem is trivial. The case
for $p=0$ follows from Dold-Thom theorem (\cite{Dold-Thom}). The
case for $p=n-1$ is due to Friedlander \cite{Friedlander1}. The case
for $p=n$ is from the definition. In particular, these invariants
are trivial for smooth projective varieties with dimension less than
or equal to two.}

\medskip
Another result is this paper is the following:

{\Th (Lawson homology for a blowup) Let $X$ be smooth projective
manifold and $Y\subset X$ be a smooth subvariety of codimension r.
Let $\sigma:\tilde{X}_Y\rightarrow X$ be the blowup of $X$ along
$Y$, $\pi:D=\sigma^{-1}(Y)\rightarrow Y$ the natural map, and
$i:D=\sigma^{-1}(Y)\rightarrow \tilde{X}_Y$ the exceptional divisor
of the blowing up. Then for each $p$, $k$ with $k\geq 2p\geq 0$, we
have the following isomorphism

  $$
\begin{array}{cc}
  &I_{p,k}: \bigg\{\bigoplus_{1\leq j \leq
r-1}L_{p-j}H_{k-2j}(Y) \bigg\}\bigoplus L_pH_k(X)\cong L_pH_k(\tilde{X}_Y)\\

\end{array}
  $$
 }

\medskip
As applications, we have the following

 {\Cor For each $n\geq 5$, there exists  a \emph{rational} manifold $X$ with
$\dim(X)=n$ such that
 $$\dim_{\mathbb{Q}}\{ {\rm
Griff}_p(X)\otimes {\mathbb{Q}}\}=\infty, \quad 2\leq p\leq n-3.$$}

{\Cor  For any integer $p>1$ and $k\geq 0$, there exists
\emph{rational} projective manifold $X$ such that
$L_pH_{k+2p}(X)\otimes {\mathbb{Q}}$ is an infinite dimensional
vector space over $\mathbb{Q}$. }

\medskip
The main tools used to prove the main result are: the long exact
localization sequence given by Lima-Filho in \cite{Lima-Filho}, the
explicit formula for the Lawson homology of codimension-one cycles
on a smooth projective manifold given by Friedlander in
\cite{Friedlander1}, and the weak factorization theorem proved by
Wlodarczyk and others in \cite{Wlodarczyk} and in \cite{AKMW}.

 \s{Some fundamental
materials in Lawson homology} {\hskip .2 in} First recall that for a
morphism $f:U\rightarrow V$ between  projective varieties, there
exist induced homomorphism
$$f_*: L_pH_k(U)\rightarrow L_pH_k(V)$$
for all $k\geq 2p\geq 0$, and if $g:V\rightarrow W$ is another
morphism between projective varieties, then $$(g\circ f)_*=g_*\circ
f_*.$$

Furthermore, it has been shown by C. Peters [Pe] that if $U$ and $V$
are smooth and projective, there are Gysin ``wrong way"
homomorphisms $f^*: L_pH_k(V)\rightarrow L_{p-c}H_{k-2c}(U)$, where
$c=\dim(V)-\dim(U)$. If $g:V\rightarrow W$ is another morphism
between smooth projective varieties, then
$$(g\circ f)^*=f^*\circ g^*.
$$

\medskip
Recall also the fact that there is a long exact sequence (cf.
\cite{Lima-Filho}, also \cite{Friedlander-Gabber})
$$\cdots\rightarrow L_pH_k(U-V)\rightarrow L_pH_k(U)\rightarrow L_pH_k(V)\rightarrow
L_pH_{k-1}(U-V)\rightarrow\cdots,$$ where $U$ is quasi-projective
and $U-V$ is any algebraic closed subset in $U$.

\medskip
Let $X$ be a smooth projective variety and $i_0:Y\hookrightarrow X$
a smooth subvariety of codimension $r\geq 2$. Let
$\sigma:\tilde{X}_Y\rightarrow X$ be the blowup of $X$ along $Y$,
$\pi:D=\sigma^{-1}(Y)\rightarrow Y$ the natural map, and
$i:D=\sigma^{-1}(Y)\hookrightarrow \tilde{X}_Y$ the exceptional
divisor of the blowup. Set $U:= X-Y\cong \tilde{X}_Y - D$. Denote by
$j_0$ the inclusion $U\subset X$ and $j$ the inclusion $U\subset
\tilde{X}_Y$. Note that $\pi:D=\sigma^{-1}(Y)\rightarrow Y$ makes
$D$ into a projective bundle of rank $r-1$, given precisely by
$D=\P(N_{Y/X})$ and we have (cf. [\cite{Voision}, pg. 271])

$$ {\mathcal{O}}_{\tilde{X}_Y}(D)|_{D}=
{\mathcal{O}}_{\P(N_{Y/X})}(-1).
$$

Denote by $h$ the class of ${\mathcal{O}}_{\P(N_{Y/X})}(-1)$ in
${\rm Pic}(D)$. We have $h=-D_{|D}$ and
$-h=i^*i_*:L_qH_m(D)\rightarrow L_{q-1}H_{m-2}(D)$ for $0\leq 2q\leq
m$ (\cite{Friedlander-Gabber}, Theorem 2.4], [\cite{Peters}, Lemma
11]). The last equality can be equivalently regarded as a Lefschetz
operator
\begin{eqnarray}
-h=i^*i_*:L_qH_m(D)\rightarrow L_{q-1}H_{m-2}(D), \quad 0\leq 2q\leq
m.
\end{eqnarray}

\medskip
The proof of the main result is based on the following lemmas:

{\Lemma For each $p\geq 0$, we have the following commutative
diagram
$$
\begin{array}{ccccccccccc}
\cdots\rightarrow & L_pH_k(D) & \stackrel{i_*}{\rightarrow} &
L_pH_k({\tilde{X}_Y})
 & \stackrel{j^*}{\rightarrow} & L_pH_k(U) & \stackrel{\delta_*}{\rightarrow} & L_pH_{k-1}(D) & \rightarrow & \cdots & \\

 & \downarrow \pi_*&   & \downarrow \sigma_*&   & \downarrow  \cong&   & \downarrow \pi_*&   &  &\\

 \cdots\rightarrow & L_pH_k(Y) & \stackrel{(i_0)*}{\rightarrow} &
L_pH_k({X}) & \stackrel{j_0^*}{\rightarrow}&
L_pH_k(U)&\stackrel{(\delta_0)_*}{\rightarrow} & L_pH_{k-1}(Y) &
\rightarrow & \cdots&

\end{array}
$$
}

 \bp This is from the corresponding commutative diagram of
 fibration sequences of p-cycles. More precisely, to show the first
 square, we begin from the following commutative diagram
$$
\begin{array}{cccc}
 D & \stackrel{i}{\hookrightarrow} & {\tilde{X}_Y}& \\
 \downarrow \pi &   & \downarrow \sigma&   \\
 Y & \stackrel{i_0}{\hookrightarrow} & {X} .&
\end{array}
$$
From this, we obtain the corresponding commutative diagram of
p-cycles:
$$
\begin{array}{cccc}
{\cal Z}_p(D) & \stackrel{i_*}{\hookrightarrow} & {\cal Z}_p({\tilde{X}_Y})& \\
\downarrow \pi_* &   & \downarrow \sigma_*&   \\
{\cal Z}_p(Y) & \stackrel{(i_0)_*}{\hookrightarrow} & {\cal
Z}_p({X}) .&
\end{array}
$$

Since $Y$ is a smooth projective variety, ${\tilde{X}_Y}$ and $D$
are smooth projective varieties, we have the following commutative
diagram

$$
\begin{array}{cccc}
  {\cal Z}_p({\tilde{X}_Y})&\rightarrow & {\cal Z}_p({\tilde{X}_Y})/{\cal Z}_p(D) \\
\downarrow \sigma_* &   & \downarrow \cong&   \\
{\cal Z}_p(X) & \rightarrow & {\cal Z}_p({X})/{\cal Z}_p(Y) .&
\end{array}
$$

Therefore we obtain the following commutative diagram of the
fibration sequences of p-cycles
$$
\begin{array}{cccccc}
 {\cal Z}_p(D) & \stackrel{i_*}{\hookrightarrow} & {\cal
 Z}_p({\tilde{X}_Y})&\rightarrow & {\cal Z}_p({\tilde{X}_Y})/{\cal Z}_p(D) \\

 \downarrow \pi_* &   & \downarrow \sigma_*&   & \downarrow \cong&   \\

 {\cal Z}_p(Y) & \stackrel{(i_0)_*}{\hookrightarrow} &  {\cal
Z}_p(X) & \rightarrow & {\cal Z}_p({X})/{\cal Z}_p(Y) .&
\end{array}
$$
where the fact that the rows are fibration sequences is due to Lima-
Filho \cite{Lima-Filho}.

 By taking the homotopy groups of these fibration sequences, we
get the long exact sequences of commutative diagram given in the
Lemma.

 \qe

 {\Prop If $p=0$, then we have the following commutative diagram
$$
\begin{array}{ccccccccccc}
\cdots\rightarrow & H_k(D) & \stackrel{i_*}{\rightarrow} &
H_k({\tilde{X}_Y})
 & \stackrel{j^*}{\rightarrow} & H_k^{BM}(U) & \stackrel{\delta_*}{\rightarrow} & H_{k-1}(D) & \rightarrow & \cdots & \\

 & \downarrow \pi_*&   & \downarrow \sigma_*&   & \downarrow  \cong&   & \downarrow \pi_*&   &  &\\

 \cdots\rightarrow & H_k(Y) & \stackrel{(i_0)*}{\rightarrow} &
H_k({X}) & \stackrel{j_0^*}{\rightarrow} &
H_k^{BM}(U)&\stackrel{(\delta_0)_*}{\rightarrow} & H_{k-1}(Y) &
\rightarrow & \cdots&

\end{array}
$$
Moreover, if $x\in H_k(D)$ maps to zero under $\pi_*$ and $i_*$,
then $x=0\in H_k(D)$. }

\medskip
 \bp The first conclusion follows directly from Lemma 2.1 with $p=0$ and the Dold-Thom Theorem. For the
 second conclusion assume $i_*(x)=0$ and $\pi_*(x)=0$. Then there exists an element $y\in
 H_{k+1}^{BM}(U)$ such that the image of $y$ under the boundary map $(\delta_0)_*:H_{k+1}^{BM}(U)\rightarrow
 H_{k}(Y)$ is 0 by the given condition. Hence there exists an
 element $z\in H_{k+1}({X})$ such that $(j_0)^*(z)=y$. Now the
 surjectivity of the map $\sigma_*:H_{k+1}({\tilde{X}_Y})\rightarrow H_{k+1}({X})$
 implies that there is an element $\tilde{z}\in
 H_{k+1}({\tilde{X}_Y})$ such that $j^*(\tilde{z})=y$. Therefore,
 $x=0 \in  H_k(D)$.

 \qe

{\Cor If $p=n-2$, then we have the following commutative diagram
$$
\begin{array}{ccccccccccc}
\cdots\rightarrow & L_{n-2}H_k(D) & \stackrel{i_*}{\rightarrow} &
L_{n-2}H_k({\tilde{X}_Y})
 & \stackrel{j^*}{\rightarrow} & L_{n-2}H_k(U) & \stackrel{\delta_*}{\rightarrow} & L_{n-2}H_{k-1}(D) & \rightarrow & \cdots & \\

 & \downarrow \pi_*&   & \downarrow \sigma_*&   & \downarrow  \cong&   & \downarrow \pi_*&   &  &\\

 \cdots\rightarrow & L_{n-2}H_k(Y) & \stackrel{(i_0)*}{\rightarrow} &
L_{n-2}H_k({X}) & \stackrel{j_0^*}{\rightarrow} &
L_{n-2}H_k(U)&\stackrel{(\delta_0)_*}{\rightarrow} &
L_{n-2}H_{k-1}(Y) & \rightarrow & \cdots&

\end{array}
$$
}

{\Lemma For each $p$, we have the following commutative diagram
$$
\begin{array}{ccccccccccc}
\cdots\rightarrow & L_pH_k(D) & \stackrel{i_*}{\rightarrow} &
L_pH_k({\tilde{X}_Y})
 & \stackrel{j^*}{\rightarrow} & L_pH_k(U) & \stackrel{\delta_*}{\rightarrow} & L_pH_{k-1}(D) & \rightarrow & \cdots & \\

 & \downarrow \Phi_{p,k}&   & \downarrow \Phi_{p,k}&   & \downarrow  \Phi_{p,k}&   & \downarrow \Phi_{p,k-1}&   &  &\\

 \cdots\rightarrow & H_k(D) & \stackrel{i_*}{\rightarrow} &
H_k({\tilde{X}_Y}) & \stackrel{j^*}{\rightarrow} &
H_k^{BM}(U)&\stackrel{\delta_*}{\rightarrow} & H_{k-1}(D) &
\rightarrow & \cdots&

\end{array}
$$
In particular, it is true for $p=1, n-2$.}

\medskip
\bp See \cite{Lima-Filho} and also \cite{Friedlander-Mazur}.
 \qe

{\Lemma For each $p$, we have the following commutative diagram
$$
\begin{array}{ccccccccccc}
\cdots\rightarrow & L_pH_k(Y) & \stackrel{(i_0)_*}{\rightarrow} &
L_pH_k({{X}})
 & \stackrel{j^*}{\rightarrow}& L_pH_k(U) & \stackrel{(\delta_0)_*}{\rightarrow} & L_pH_{k-1}(Y) & \rightarrow & \cdots & \\

 & \downarrow \Phi_{p,k}&   & \downarrow \Phi_{p,k}&   & \downarrow  \Phi_{p,k}&   & \downarrow \Phi_{p,k-1}&   &  &\\

 \cdots\rightarrow & H_k(Y) & \stackrel{(i_0)_*}{\rightarrow} &
H_k({X}) & \stackrel{j^*}{\rightarrow} &
H_k^{BM}(U)&\stackrel{(\delta_0)_*}{\rightarrow} & H_{k-1}(Y) &
\rightarrow & \cdots&

\end{array}
$$
In particular, it is true for $p=1, n-2$.}
\\

\bp See \cite{Lima-Filho} and also \cite{Friedlander-Mazur}.
\qe


\s{Lawson homology for blowups}

\hskip .2in
 As an application of Lemma 2.1, we give an explicit formula for a
 blowup in Lawson homology. Since it may have some independent interest, we devote a separate section
 to it. First, we want to revise the projective bundle theorem given
 by Friedlander and Gabber (\cite{Friedlander-Gabber}, Prop.2.5). It is convenient to extend the
 definition of Lawson homology by setting

 $$ L_pH_k(X)=L_0H_k(X), \quad if \quad p<0.$$

 Now we have the following revised ``Projective Bundle Theorem":

 {\Prop Let $E$ be an algebraic vector bundle of rank $r$ over a
 smooth projective variety
 $Y$, then for each $p\geq 0$ we have
 $$ L_pH_k(\P(E))\cong \bigoplus_{j=0}^{r-1}L_{p-j}H_{k-2j}(Y)$$
 where $\P(E)$ is the projectivization of the vector bundle $E$.
 }

{\Rem The difference between this and the projective bundle theorem
of \cite{Friedlander-Gabber} is that here we place no restriction on
$p$.}

\medskip
\bp For $p\geq r-1$, this is exactly the projective bundle theorem
given in \cite{Friedlander-Gabber}. If $p<r-1$, we have the same
method of \cite{Friedlander-Gabber}, i.e., the localization sequence
and the naturality of $\Phi$, to reduce to the case in which $E$ is
trivial. From
$${\cal{Z}}_0(\P^{r-1}\times Y)\rightarrow {\cal{Z}}_0(\P^r\times
Y)\rightarrow {\cal{Z}}_0(\mathbb{C}^r\times Y),$$ we have the long
exact localization sequence given at the beginning of section 2:
$$\cdots\rightarrow L_0H_k(\P^{r-1}\times Y)\rightarrow
L_0H_k(\P^r\times Y)\rightarrow L_0H_k(\mathbb{C}^r\times
Y)\rightarrow
 L_0H_{k-1}(\P^{r-1}\times Y)\rightarrow\cdots.
$$

From this, and the K\"unneth formula for $\P^r\times Y$, we have the
following isomorphism:
$$(*)\quad H_{k-2r}(Y)\cong
L_0H_{k}(\mathbb{C}^r\times Y)\cong H_{k}^{BM}(\mathbb{C}^r\times
Y).
$$

Note that
$$(**)\quad\quad H_{k-2r}(Y)\cong L_{p-r}H_{k-2r}(Y) \quad {\rm if}\quad p\leq r.
$$

All the remaining arguments are the same as those in
[\cite{Friedlander-Gabber}, Prop 2.5], as we review in the
following.

We want to use induction on $r$. For $r-1=p$, the conclusion holds.
From the commutative diagram of abelian groups of cycles:
$$
\begin{array}{cccc}
\{\oplus_{j=0}^{p}Z_{p-j}(X)\}\bigoplus\{\oplus_{j=p+1}^{r-1}
Z_0(X\times \mathbb{C}^{j-p})\}&\rightarrow &
\{\oplus_{j=0}^{p}Z_{p-j}(X)\}\bigoplus\{\oplus_{j=p+1}^{r}
Z_0(X\times \mathbb{C}^{j-p})\}\\
\downarrow & & \downarrow  \\
 Z_{p}(X\times \P^{r-1})           &\rightarrow &   Z_{p}(X\times \P^{r})
\end{array}
$$

We obtain the commutative diagram of fibration sequences:
$$
\begin{array}{cccccc}
\{\oplus_{j=0}^{p}Z_{p-j}(X)\}\bigoplus\{\oplus_{j=p+1}^{r-1}Z_{p-j}(X)\}&\rightarrow&
\{\oplus_{j=0}^{p}Z_{p-j}(X)\}\bigoplus\{\oplus_{j=p+1}^{r}Z_{p-j}(X)\}&\rightarrow&
Z_0(X\times \mathbb{C}^{r-p})\\

\downarrow & & \downarrow &&\downarrow  \\
 Z_{p}(X\times \P^{r-1})  &\rightarrow &   Z_{p}(X\times
 \P^{r})&\rightarrow& Z_p(X\times \mathbb{C}^{r})
\end{array}
$$
where $Z_{p-j}(X):=Z_0(X\times \mathbb{C}^{j-p})$ for $p-j<0$.

The first vertical arrow is a homotopy equivalence by induction. The
last one is a homotopy equivalence by Complex Suspension Theorem
\cite{Lawson1}. Hence by the Five Lemma, we obtain the homotopy
equivalence of the middle one.

The proof is completed by combining this with (*) and (**) above.

\qe

{\Rem The isomorphism
$$ \psi: \bigoplus_{j=0}^{r-1}L_{p-j}H_{k-2j}(Y)\stackrel{\cong}{\longrightarrow} L_pH_k(\P(E))$$
in Proposition 3.1 is given explicitly by
$$ \psi(u_0,u_1,\cdots,u_{r-1})=\sum_{j=0}^{r-1}h^j\pi^*u_j
$$
where $h$ is the Lefschetz hyperplane operator $h:
L_qH_m(\P(E))\rightarrow L_{q-1}H_{m-2}(\P(E))$ defined in (1). For
$p\geq r-1$, this explicit formula has been proved in
[\cite{Friedlander-Gabber}, Prop. 2.5]. In the remaining cases, $h$
is the Lefschetz hyperplane operator $h: H_m(\P(E))\rightarrow
H_{m-2}(\P(E))$ defined in (1).

}

\medskip
Using the notations in section 2, we have the following:

{\Th (Lawson homology for a blowup) Let $X$ be smooth projective
manifold and $Y\subset X$ be a smooth subvariety of codimension r.
Let $\sigma:\tilde{X}_Y\rightarrow X$ be the blowup of $X$ along
$Y$, $\pi:D=\sigma^{-1}(Y)\rightarrow Y$ the natural map, and
$i:D=\sigma^{-1}(Y)\rightarrow \tilde{X}_Y$ the exceptional divisor
of the blowing up. Then for each $p$, $k$ with $k\geq 2p\geq 0$, we
have the following isomorphism

$$
I_{p,k}: \bigg\{\bigoplus_{1\leq j \leq r-1}L_{p-j}H_{k-2j}(Y)
\bigg\}\oplus L_pH_k(X)\stackrel{\cong}{\longrightarrow}
L_pH_k(\tilde{X}_Y)
$$
given by

$$
I_{p,k}(u_1,\cdots,u_{r-1},
u)=\sum_{j=1}^{r-1}i_*h^j\pi^*u_j+\sigma^*u
$$
 }

\bp Let $U:=\tilde{X}_Y-D=X-Y$. By the definitions of the maps $i$,
$\pi$ and $\sigma$, and Lemma 2.1, we have  the following
commutative diagram of the long exact localization sequences:
\begin{eqnarray}
\begin{array}{ccccccccccc}
\cdots\rightarrow & L_pH_k(D) & \stackrel{i_*}{\rightarrow} &
L_pH_k({\tilde{X}_Y})
 & \stackrel{j^*}{\rightarrow} & L_pH_k(U) & \stackrel{\delta_*}{\rightarrow} & L_pH_{k-1}(D) & \rightarrow & \cdots & \\

 & \downarrow \pi_*&   & \downarrow \sigma_*&   & \downarrow  \cong&   & \downarrow \pi_*&   &  &\\

 \cdots\rightarrow & L_pH_k(Y) & \stackrel{(i_0)*}{\rightarrow} &
L_pH_k({X}) & \stackrel{j_0^*}{\rightarrow}&
L_pH_k(U)&\stackrel{(\delta_0)_*}{\rightarrow} & L_pH_{k-1}(Y) &
\rightarrow & \cdots&
\end{array}
\end{eqnarray}

From this, and the surjectivity of $j^*$, we have

$$ L_pH_{2p}(\tilde{X}_Y)=\sigma^*L_pH_{2p}(X)+i_*L_pH_{2p}(D).
$$

By the ``revised" projective bundle theorem above, for any $p\geq
0$, there is an isomorphism
 $$ L_pH_k(D)\cong\bigoplus_{j=0}^{r-1} h^j\pi^* L_{p-j}H_{k-2j}(Y),\quad 0\leq 2p\leq k.$$

Hence we see that
\begin{eqnarray}
L_pH_{2p}(\tilde{X}_Y)=\sigma^*L_pH_{2p}(X)+ \Sigma_{j=0}^{r-1}
i_*h^j\pi^* L_{p-j}H_{2p-2j}(Y).
\end{eqnarray}

But clearly by Lemma 2.1 and the projective bundle theorem, if $u\in
L_pH_{k}(Y)$, then
\begin{eqnarray*}
\sigma_*(i_*h^{r-1}\pi^*(u))=(i_0)_*(u).
\end{eqnarray*}

Since $\sigma$ is a birational morphism, it has degree one.  As a
directly corollary of the projection formula (cf. \cite{Peters},
Lemma 11 c.), we have $\sigma_*(\sigma^* a)=a$ for any $a\in
L_pH_k(X)$. We have
$$\sigma_*(\sigma^* ((i_0)_*u))=(i_0)_*u, \quad u\in L_pH_{k}(Y).
$$

Thus we obtain the relations
$$v:=i_* h^{r-1}\pi^*u-\sigma^* ((i_0)_*u)\in \ker \sigma_*, \quad u\in L_pH_{k}(Y)
$$

Since $j^*=(j_0)^*\sigma_*$ in (2), we get $j^*(v)=0$. From the
exactness of the upper row in (2), we get
\begin{eqnarray}
v\in \sum_{j=1}^{r-1}i_*h^jL_{p-j}H_{k-2j}(Y).
\end{eqnarray}

The equality (3) and the relation (4) together imply immediately
that the map $I_{p,2p}$ is surjective for the case $k=2p$.

To prove the injectivity for the case that $k=2p$, consider $(u_1,
u_2,\cdots,u_{r-1}, u)\in\ker I_{p,2p}$. Applying $\sigma_*$, we
find that $u=0$. Note that $i^*i_*=-h$. Applying $i^*$ to the
equality
$$
\sum_{j=1}^{r-1}i_*h^j\pi^*u_j=0,
$$
we get
$$
\sum_{j=1}^{r-1}h^{j+1}\pi^*u_j=0 \in L_{p-1}H_{k-2}(D).
$$

The isomorphism in Proposition 3.1 implies that $u_j=0$ for $1\leq
j\leq r-1$. This completes the proof for the case $k=2p$.

From this and  (2), we have
\begin{eqnarray}
\begin{array}{ccccccccccc}
\cdots\rightarrow & L_pH_{2p+1}(D) & \stackrel{i_*}{\rightarrow} &
L_pH_{2p+1}({\tilde{X}_Y})
 & \stackrel{j^*}{\rightarrow} & L_pH_{2p+1}(U) & \stackrel{\delta_*}{\rightarrow} & 0 & \  &  & \\

 & \downarrow \pi_*&   & \downarrow \sigma_*&   & \downarrow  \cong&   &      &  &\\

 \cdots\rightarrow & L_pH_{2p+1}(Y) & \stackrel{(i_0)*}{\rightarrow} &
L_pH_{2p+1}({X}) & \stackrel{j_0^*}{\rightarrow}&
L_pH_{2p+1}(U)&\stackrel{(\delta_0)_*}{\rightarrow} & 0 &  & &
\end{array}
\end{eqnarray}

Now the situation for $k=2p+1$ is the same as that in the case
$k=2p$. From $(5)$ and the ``revised" projective bundle theorem, we
have

\begin{eqnarray}
L_pH_{2p+1}(\tilde{X}_Y)=\sigma^*L_pH_{2p+1}(X)+ \Sigma_{j=0}^{r-1}
i_*h^j\pi^* L_{p-j}H_{2p+1-2j}(Y).
\end{eqnarray}

From (4) and (6), we obtain the surjectivity of $I_{p,2p+1}$ for the
case that $k=2p+1$.

To prove the injectivity, consider $(u_1, u_2,\cdots,u_{r-1},
u)\in\ker I_{p,2p+1}$. Applying $\sigma_*$, we find that $u=0$. Note
that $i^*i_*=-h$. By applying $i^*$ to the equality
$$
\sum_{j=1}^{r-1}i_*h^j\pi^*u_j=0,
$$
we get
$$
\sum_{j=1}^{r-1}h^{j+1}\pi^*u_j=0 \in L_{p-1}H_{k-2}(D).
$$

The isomorphism in Proposition 3.1 again implies that $u_j=0$ for
$1\leq j\leq r-1$. This completes the proof for the case $k=2p+1$.

Now for $k\geq 2p+2$, we reach the same situation as those in the
case that $k=2p$ or $k=2p+1$. More precisely, we give the complete
argument by using mathematical induction.

Suppose that we have \begin{eqnarray}
\begin{array}{ccccccccccc}
\cdots\rightarrow & L_pH_{2p+m}(D) & \stackrel{i_*}{\rightarrow} &
L_pH_{2p+m}({\tilde{X}_Y})
 & \stackrel{j^*}{\rightarrow} & L_pH_{2p+m}(U) & \stackrel{\delta_*}{\rightarrow} & 0 & \  &  & \\

 & \downarrow \pi_*&   & \downarrow \sigma_*&   & \downarrow  \cong&   &      &  &\\

 \cdots\rightarrow & L_pH_{2p+m}(Y) & \stackrel{(i_0)*}{\rightarrow} &
L_pH_{2p+m}({X}) & \stackrel{j_0^*}{\rightarrow}&
L_pH_{2p+m}(U)&\stackrel{(\delta_0)_*}{\rightarrow} & 0 &  & &
\end{array}
\end{eqnarray}
for some integer $m\geq 0$.

We want to prove that $I_{p,2p+m}$ is an isomorphism and
\begin{eqnarray}
\begin{array}{ccccccccccc}
\cdots\rightarrow & L_pH_{2p+m+1}(D) & \stackrel{i_*}{\rightarrow} &
L_pH_{2p+m+1}({\tilde{X}_Y})
 & \stackrel{j^*}{\rightarrow} & L_pH_{2p+m+1}(U) & \stackrel{\delta_*}{\rightarrow} & 0 & \  &  & \\

 & \downarrow \pi_*&   & \downarrow \sigma_*&   & \downarrow  \cong&   &      &  &\\

 \cdots\rightarrow & L_pH_{2p+m+1}(Y) & \stackrel{(i_0)*}{\rightarrow} &
L_pH_{2p+m+1}({X}) & \stackrel{j_0^*}{\rightarrow}&
L_pH_{2p+m+1}(U)&\stackrel{(\delta_0)_*}{\rightarrow} & 0 &  & &
\end{array}
\end{eqnarray}

Once this step is done, it completes the proof of the theorem.

From the assumption (7), we have
\begin{eqnarray}
L_pH_{2p+m}(\tilde{X}_Y)=\sigma^*L_pH_{2p+m}(X)+ \Sigma_{j=0}^{r-1}
i_*h^j\pi^* L_{p-j}H_{2p+m-2j}(Y).
\end{eqnarray}

From (4) for $k=2p+m$ and (9), we obtain the surjectivity of
$I_{p,2p+m}$ for the case that $k=2p+m$.

To prove the injectivity, consider $(u_1, u_2,\cdots,u_{r-1},
u)\in\ker I_{p,2p+m}$. Applying $\sigma_*$, we find that $u=0$. Note
that $i^*i_*=-h$. By applying $i^*$ to the equality
$$
\sum_{j=1}^{r-1}i_*h^j\pi^*u_j=0,
$$
we get
$$
\sum_{j=1}^{r-1}h^{j+1}\pi^*u_j=0 \in L_{p-1}H_{k-2}(D).
$$

The isomorphism in Proposition 3.1 once again implies that $u_j=0$
for $1\leq j\leq r-1$. This completes the proof for the case
$k=2p+m$. Now (7) automatically reduces to (8) and  this completes
the proof of the theorem. \qe



As an application, this result gives many examples of smooth
projective manifolds (even rational ones) for which the Griffiths
group of p-cycles is infinitely generated (even modulo torsion) for
$p\geq 2$. Recall that the Griffiths group ${\rm Griff}_p(X)$ is
defined to be the p-cycles homologically equivalent to zero modulo
the subgroup of p-cycles algebraically equivalent to zero.

\medskip
\noindent{\bf Example:} Note the fact in \cite{Friedlander1} that
${\rm Griff}_2(\tilde{X}_Y)\cong L_2H_4(\tilde{X}_Y)_{hom}$.
 For $X=\P^5$,
$Y\subset \P^4$ the general hypersurface of degree 5, we obtain an
infinite dimensional ${\mathbb{Q}}$-vector space ${\rm
Griff}_2(\tilde{X}_Y)\otimes {\mathbb{Q}}$ from the fact $\dim_{
\mathbb{Q}}({\rm Griff}_1(Y)\otimes {\mathbb{Q}})=\infty$ (cf.
\cite{Clemens}). It gives the example mentioned in Remark 1.1.

From the blowup formula for Lawson homology and Clemens' result
\cite{Clemens}, we have the following
 {\Cor For each $n\geq 5$, there exists a rational manifold $X$ with $\dim(X)=n$ such that
$$\dim_{\mathbb{Q}} {\rm
Griff}_p(X)\otimes {\mathbb{Q}}=\infty, \quad 2\leq p\leq n-3.$$}

\bp Note that ${\rm Griff}_p(X)\cong L_pH_{2p}(X)_{hom}$ for any
smooth projective variety $X$. Now the remaining argument is the
direct result of Theorem 3.1 and  the above result of Clemens
\cite{Clemens}.

\qe

More generally, from the blowup formula for Lawson homology and a
result given by the author \cite{author}, we have the following

{\Cor  For any integers $p>1$ and $k\geq 0$, there exists a
\emph{rational} projective manifold $X$ such that
$L_pH_{k+2p}(X)\otimes {\mathbb{Q}}$ is infinite dimensional vector
space over $\mathbb{Q}$. }

\medskip
\bp It follows from the blowup formula for Lawson homology and
Theorem 1.4 in \cite{author}. For example, if $p=2$, $k=1$, we can
find a rational projective manifold $X$ with $\dim(X)=6$ such that
$L_2H_5(X) \otimes {\mathbb{Q}}$ is infinite dimensional
$\mathbb{Q}$-vector space.

\qe

\s{The proof of the main theorem}

\hskip .2in The following result will be used  several times in the
proof of our main theorem:

{\Th (Friedlander \cite{Friedlander1}) Let $X$ be any smooth
projective variety of dimension $n$. Then we have the following
isomorphisms
$$
\left\{
\begin{array}{l}
 L_{n-1}H_{2n}(X)\cong \Z,\\
 L_{n-1}H_{2n-1}(X)\cong H_{2n-1}(X,\Z),\\
 L_{n-1}H_{2n-2}(X)\cong H_{n-1,n-1}(X,\Z)=NS(X)\\
L_{n-1}H_{k}(X)=0 \quad for\quad k> 2n.\\

\end{array}
\right.
$$}
\qe

{\Rem In the following, we adopt the notational convention
$H_k(X)=H_k(X,\Z)$.}
\\

Now we begin the proof of our main results. There are two parts of
the proof of the main theorem: $p=1$ and $p=n-2$.\\

\noindent\textbf{Proof of the main theorem ($p=1$):}

\medskip
\noindent \textbf{Case A:}
$\sigma_*:L_{1}H_{k}(\tilde{X}_Y)_{hom}\rightarrow
L_{1}H_{k}(X)_{hom}$ is injective.

\medskip
We will use the commutative diagrams in Lemma 2.1--2.3.

Let $a\in L_{1}H_{k}({\tilde{X}_Y})_{hom}$ be such that
$\sigma_*(a)=0$. By Lemma 2.1, we have $j^*(a)=0\in L_{1}H_{k}(U)$
and hence there exists an element $b\in L_{1}H_{k}(D)$ such that
$i_*(b)=a$. Set $\tilde{b}=\pi_*(b)$. By the commutative diagram in
Lemma 2.1 again, we have $(i_0)_*(\tilde{b})=0\in L_{1}H_{k}(X)$. By
the exactness of the rows in the commutative diagram, there exists
an element $\tilde{c}\in L_1H_{k+1}(U)$ such that the image of
$\tilde{c}$ under the boundary map $(\delta_0)_*:
L_{1}H_{k+1}(U)\rightarrow L_{1}H_{k}(Y)$ is $\tilde{b}$. Note that
$\delta_*$  is the other boundary map $\delta_*:
L_{1}H_{k+1}(U)\rightarrow L_{1}H_{k}(D)$. Therefore,
$\pi_*(b-\delta_*(\tilde{c}))=0\in L_{1}H_{k}(Y)$ and
$j_*(b-\delta_*(\tilde{c}))=a$. Now by the ``revised" Projective
Bundle Theorem and Dold-Thom theorem ([DT]), we have
$L_{1}H_{k}(D)\cong L_{1}H_{k}(Y)\oplus L_{0}H_{k-2}(Y)\oplus
H_{k-4}(Y)\oplus\cdots\cong L_{1}H_{k}(Y)\oplus H_{k-2}(Y)\oplus
H_{k-4}(Y)\oplus\cdots$. We know $b-\delta_*(\tilde{c})\in
H_{k-2}(Y)\oplus H_{k-4}(Y)\oplus\cdots$. By the explicit formula of
the cohomology (and homology) for a blowup ([GH]), we know each map
$H_{k-2*}(Y)\rightarrow H_{k}(\tilde{X}_Y)$ is injective. Hence $a$
must be zero in $L_{1}H_{k}({\tilde{X}_Y})$. This is the injectivity
of $\sigma_*$.

\medskip
\noindent\textbf{Case B:}
$\sigma_*:L_{1}H_{k}(\tilde{X}_Y)_{hom}\rightarrow
L_{1}H_{k}(X)_{hom}$ is surjective.

\medskip
Let $a\in L_{1}H_{k}(X)_{hom}$. From the surjectivity of the map
$\sigma_*:L_{1}H_{k}(\tilde{X}_Y)\rightarrow L_{1}H_{k}(X)$, there
exists an element $\tilde{a}\in L_{1}H_{k}(\tilde{X}_Y)$ such that
$\sigma_*(\tilde{a})=a$. Set $\tilde{b}=\Phi_{1,k}(\tilde{a})$. By
the commutative diagram in Lemma 2.1, we have $j^*(\tilde{b})=0\in
H_k^{BM}(U)$. From the exactness of the rows of the diagram in Lemma
2.1, we have an element $\tilde{c}\in H_{k}(D)$ such that
$i_*(\tilde{c})=\tilde{b}$. Set $c=\pi_*(\tilde{c})$. Then
$(i_0)_*(c)=0$ by the assumption of $a$ and the commutative of the
diagram in Lemma 2.1. Using the exactness of rows in Lemma 2.1
again, we can find an element $d\in H_{k+1}^{BM}(U)$ such that
$(\delta_0)_*(d)=c$. Hence $i_*(\tilde{c}-\delta_*(d))=\tilde{b}\in
H_{k}(\tilde{X}_Y)$ and $\pi_*(\tilde{c}-\delta_*(d))=0$. Now we
need to use the formula $L_{1}H_{k}(D)\cong  L_{1}H_{k}(Y)\oplus
H_{k-2}(Y)\oplus H_{k-4}(Y)\oplus\cdots$ again. From this we can
find an element $e\in L_{1}H_{k}(D)$ such that
$\Phi_{1,k}(e)=\tilde{c}-\delta(d)$. Obviously,
$\Phi_{1,k}(\tilde{a}-i_*(e))=0$ and $\sigma_*(\tilde{a}-i_*(e))=a$
as we want.

 \qe

\noindent\textbf{Proof of the main theorem ($p=n-2$):}

\medskip
\noindent\textbf{Case 1:} $\sigma_*$ is injective.

\medskip
The injectivity of $j_0^*:L_{n-2}H_{k}({X})_{hom}\rightarrow
L_{n-2}H_{k}(U)_{hom}$ is trivial since the dim$(Y)\leq n-2$, where
$j_0:U\rightarrow X$ is the inclusion. In fact, if dim$(Y)<n-2$,
$j_0^*:L_{n-2}H_{k}({X})\rightarrow L_{n-2}H_{k}(U)$ is an
isomorphism and so is $j_0^*:L_{n-2}H_{k}({X})_{hom}\rightarrow
L_{n-2}H_{k}(U)_{hom}$. If dim$(Y)=n-2$, then for $k\geq 2(n-2)+1$
the injectivity of $j_0^*$ is from the commutative diagram in Lemma
2.2, and the vanishing of $ L_{n-2}H_{k}(Y)$ and $H_{k}(Y)$; for
$k=2(n-2)$, the injectivity of $j_0^*$ is from the commutative
diagram in Lemma 2.2, and the nontriviality of $(i_0)_*:
H_{2(n-2)}(Y)\rightarrow H_{2(n-2)}({X})$, since $Y$ is a K\"{a}hler
submanifold of $X$ with complex dimension $n-2$.

Now we need to prove $j^*:L_{n-2}H_{k}({\tilde{X}_Y})_{hom}\rightarrow
L_{n-2}H_{k}(U)_{hom}$ is injective, where $j:U\rightarrow {\tilde{X}_Y}$ the
inclusion.
Let $a\in L_{n-2}H_{k}({\tilde{X}_Y})_{hom}$ such that $j^*(a)=0 \in
L_{n-2}H_{k}(U)_{hom}$, then there exists an element $b\in
L_{n-2}H_{k}(D)$ such that $i_*(b)=a$. Now by the commutative
diagram in Corollary 2.1, we have $j_0^*(\sigma_*(a))=0$. Set
$a^{\prime}\equiv \sigma_*(a)$. From the exactness of localization
sequence in the bottom row in Corollary 2.1, there is an element
$b^{\prime}\in L_{n-2}H_{k}(Y)$ such that
$(i_0)_*(b^{\prime})=a^{\prime}$.

\medskip
\begin{enumerate}
\item[]{\textbf{Claim:}}  In the commutative diagram in Corollary 2.1,
there exists an element $c^{\prime}\in L_{n-2}H_{k+1}(U)$ such that
$(\delta_0)_*(c')=b'$ under the map $(\delta_0)_*:
L_{n-2}H_{k+1}(U)\rightarrow L_{n-2}H_{k}(Y)$ and $\delta_*(c')=b$
under the map $\delta_*: L_{n-2}H_{k+1}(U)\rightarrow
L_{n-2}H_{k}(D)$.

\medskip
\noindent{\rm \textbf{Proof of the claim:}} Since $\Phi_{n-2,k}:
L_{n-2}H_{k}({Y})\cong H_{k}(Y)$ (note: $k\geq 2(n-2)\geq {\rm
dim}(Y)$), we use the same notation $b^{\prime}$ for its image in
$H_{k}(Y)$ since $L_{n-2}H_{k}(Y)\rightarrow H_{k}(Y)$ is injective
for all $k\geq 2(n-2)$. At the beginning of the proof of the
injectivity of the main theorem, we have already shown that
$j_0^*:L_{n-2}H_{k}({X})_{hom}\rightarrow L_{n-2}H_{k}(U)_{hom}$ is
injective. That is to say, $(i_0)_*(b^{\prime})=0\in
L_{n-2}H_{k}(X)_{hom}$. Hence there exists an element $c\in
L_{n-2}H_{k+1}(U)$ such that whose image is $b^{\prime}$ under the
boundary map $(\delta_0)_*: L_{n-2}H_{k+1}(U)\rightarrow
L_{n-2}H_{k}(Y)$. Let $\tilde{b}$ be the image of $c$ under the map
$L_{n-2}H_{k+1}(U)\rightarrow L_{n-2}H_{k}(D)$. Now
$\pi_*(\tilde{b}-b)=0\in L_{n-2}H_{k}(Y)$ and
$i_*(\Phi_{n-2,k}(\tilde{b}-b))=0\in H_{k}({\tilde{X}_Y})$, by
Proposition 2.1, we have  $\Phi_{n-2,k}(\tilde{b}-b)=0$. Since
$\Phi_{n-2,k}$ is injective on $L_{n-2}H_{k}(D)$ (see Theorem 4.1),
we get $\tilde{b}-b=0$. This $c$ satisfies both conditions of the
claim.

\qe
\end{enumerate}

\medskip
Now everything is clear. The element $a$ comes from the element $c$
in $L_{n-2}H_{k+1}(U)$. By the exactness of the localization
sequence in the upper row in Lemma 2.1, we get $a=0 \in
L_{n-2}H_{k}({\tilde{X}_Y})$. This completes the proof of the
injectivity.

\medskip
\noindent\textbf{Case 2:} $\sigma_*$ is surjective.

\medskip
Similar to the injectivity, the surjectivity of
$j_0^*:L_{n-2}H_{k}({X})_{hom}\rightarrow L_{n-2}H_{k}(U)_{hom}$ is
trivial since the dim$(Y)\leq n-2$, where $j_0:U\rightarrow X$ is
the inclusion. In fact, if dim$(Y)<n-2$,
$j_0^*:L_{n-2}H_{k}({X})\rightarrow L_{n-2}H_{k}(U)$ is an
isomorphism and so is $j_0^*:L_{n-2}H_{k}({X})_{hom}\rightarrow
L_{n-2}H_{k}(U)_{hom}$. If dim$(Y)=n-2$, then the surjectivity of
$j_0^*$ is from the commutative diagram in Lemma 2.3, and the
isomorphism $\Phi_{n-2,2(n-2)}: L_{n-2}H_{2(n-2)}(Y)\cong
H_{2(n-2)}(Y)\cong \Z$.

We only need to show $j^*:L_{n-2}H_{k}({\tilde{X}_Y})_{hom}\cong
L_{n-2}H_{k}(U)_{hom}$, where $j:U\rightarrow {\tilde{X}_Y}$ the
inclusion. There are a few cases.
\begin{enumerate}
\item[(a)] For the case that $k=2(n-2)$, the map
$j^*:L_{n-2}H_{k}({\tilde{X}_Y})\rightarrow L_{n-2}H_{k}(U)$ is a
surjective map. Hence the induced map
$j^*$ on $L_{n-2}H_{k}({\tilde{X}_Y})_{hom}$ is also surjective by
trivial reason.
\item[(b)] The case that $k=2(n-2)+1$. By the commutative diagram in
Lemma 2.2, and note that the map $\Phi_{n-2,2(n-2)}:
L_{n-2}H_{2(n-2)}(D)\rightarrow H_{2(n-2)}(D)$ is injective, we
have, for $a\in L_{n-2}H_{2(n-2)+1}(U)_{hom}$, the image of $a$
under the boundary map $\delta_*: L_{n-2}H_{2(n-2)+1}(U)
\rightarrow L_{n-2}H_{2n}(D)$ must be zero. Hence $a$ comes from an
element $b\in L_{n-2}H_{2(n-2)+1}({\tilde{X}_Y})$. If
$\bar{b}:=\Phi_{n-2,2(n-2)+1}(b)\neq 0$, then $\exists c\in
L_{n-2}H_{2(n-2)+1}(D)$ such that $b-i_*(c)\in
L_{n-2}H_{2(n-2)+1}({\tilde{X}_Y})_{hom}$ and $j^*(b-i_*(c))=a$. In
fact, since $j^*(\bar{b})=0$, there exists $\bar{c} \in
H_{2(n-2)+1}(D)$ such that $(i_0)_*(\bar{c})=\bar{b}$. Note that
$\Phi_{n-2,2(n-2)+1}: L_{n-2}H_{2(n-2)+1}(D)\rightarrow
H_{2(n-2)+1}(D)$ is an isomorphism by Theorem 4.1, then there exists
$c\in L_{n-2}H_{2(n-2)+1}(D)$ such that
$\Phi_{n-2,2(n-2)+1}(c)=\bar{c}$. This shows the surjectivity in
this case.

\item[(c)] Now we only need to consider the situation that $k\geq
2(n-2)+2$. In this case, the surjectivity of
$j^*:L_{n-2}H_{k}({\tilde{X}_Y})_{hom}\rightarrow
L_{n-2}H_{k}(U)_{hom}$  is from the commutative diagram in Lemma
2.2, and the surjectivity of the map $\Phi_{n-2,k}:
L_{n-2}H_k(D)\rightarrow  H_k(D)$ (see Theorem 4.1). In fact, if
$a\in L_{n-2}H_{k}(U)_{hom}$, then by the exactness of the
commutative diagram in Lemma 2.2, there is an element $b\in
L_{n-2}H_{k}({\tilde{X}_Y})$ such that $j^*(b)=a$. Set
$\bar{b}=\Phi_{n-2,k}(b)$. Since $j^*(\bar{b})=0\in H_k^{BM}(U)$,
$\exists \bar{c}\in H_k(D)$ such that $i_*(\bar{c})=\bar{b}$. Now
$\Phi_{n-2,k}:L_{n-2}H_k(D)\cong H_k(D)$ (See Theorem 4.1), there
exists $c\in L_{n-2}H_k(D)$ such that $\Phi_{n-2,k}(c)=\bar{c}$. The
commutative diagram in Lemma 2.2 implies that
$\Phi_{n-2,k}(b-i_*(c))=0$, i.e., $b-i_*(c)\in
L_{n-2}H_k({\tilde{X}_Y})_{hom}$. The exactness of the upper row in
Lemma 2.2 gives $j^*(b-i_*(c))=a$. This completes the surjectivity
in this case.
\end{enumerate}

This completes the proof for a blow-up along a smooth subvariety
$Y$ of codimension at least 2 in $X$.

Now recall the weak factorization Theorem proved in \cite{AKMW} (and
also \cite{Wlodarczyk}) as follows:

{\Th (\cite{AKMW} Theorem 0.1.1, \cite{Wlodarczyk})
 Let $\varphi \colon X \to X'$ be a birational map of smooth
complete varieties over an algebraically closed field of
characteristic zero, which is an isomorphism over an open set $U$.
Then $f$ can be factored as a sequence of birational maps
$$X = X_0
\stackrel{\varphi_1}{\rightarrow} X_1
\stackrel{\varphi_2}{\rightarrow}\dots
\stackrel{\varphi_{n+1}}{\rightarrow} X_n = X'$$ where each $X_i$ is
a smooth complete variety, and $\varphi_{i+1}: X_i \to X_{i+1}$ is
either a blowing-up or a blowing-down of a smooth subvariety
disjoint from $U$. }

\medskip
 Note that $\varphi:X\rightarrow X^{\prime}$ is birational between projective
manifolds. We complete the proof of for the birational invariance
of $L_{n-2}H_k(X)_{hom}$ for any smooth $X$ by applying the above
theorem.

 \qe
{\Rem Griffiths \cite{Griffiths} showed the nontriviality of the
Griffiths group of 1-cycles of general quintic hypersurfaces in
$\P^4$ and Friedlander \cite{Friedlander1} showed that
$L_1H_2(X)_{hom}\cong {\rm Griff_1}(X)$ for any smooth projective
variety $X$. Hence, in general, this is a {\bf nontrivial}
birational invariant even for projective threefolds.
 }

\begin{center}{\bf Acknowledge}\end {center} I
would like to express my gratitude to my advisor, Blaine Lawson,
for all his help.

\medskip

\noindent
Department of Mathematics,\\
Stony Brook University, SUNY,\\
Stony Brook, NY 11794-3651\\
Email:wenchuan@math.sunysb.edu

\end{document}